\documentclass[a4paper,12pt]{amsart}
\usepackage{amssymb}
\usepackage{ifthen}
\usepackage[usenames]{color}
\usepackage{graphicx}
\nonstopmode \numberwithin{equation}{section}
\newtheorem{thm}{Theorem}[section]

\newtheorem{cor}{Corollary}[section]
\newtheorem{lem}{Lemma}[section]
\newtheorem{prop}{Proposition}[section]

\newtheorem{claim}{Claim}[section]
\newtheorem{conj}[equation]{Conjecture}

\theoremstyle{definition}
\newtheorem{defn}{Definition}[section]
\newtheorem{examp}{Example}[section]
\newtheorem{prob}[equation]{Problem}
\newtheorem{ques}[equation]{Question}
\newtheorem{rem}{Remark}[section]

\newtheorem{case}{Case}[section]

\newtheorem*{mysolution}{Solution}

\newcounter {own}
\def\theown {\thesection       .\arabic{own}}

\newenvironment{pf}[1][]{%
 \vskip 3mm
 \noindent
 \ifthenelse{\equal{#1}{}}%
  {{\slshape Proof. }}%
  {{\slshape #1.} }%
 }%
{\qed\bigskip}

\newcounter{alphabet}

\newenvironment{Thm}[1][]{\refstepcounter{alphabet}%
\bigskip%
\noindent%
{\bf Theorem \Alph{alphabet}}%
\ifthenelse{\equal{#1}{}}{}{ (#1)}%
{\bf .} \itshape}{\vskip 8pt}

\newcommand{\IC}{{\mathbb C}}

\newcommand{\ID}{{\mathbb D}}

\def\be{\begin{equation}}
\def\ee{\end{equation}}

\newcommand{\ben}{\begin{enumerate}}
\newcommand{\een}{\end{enumerate}}

\newcommand{\blem}{\begin{lem}}
\newcommand{\elem}{\end{lem}}
\newcommand{\bthm}{\begin{thm}}
\newcommand{\ethm}{\end{thm}}
\newcommand{\bcor}{\begin{cor}}
\newcommand{\ecor}{\end{cor}}
\newcommand{\beg}{\begin{examp}}
\newcommand{\eeg}{\end{examp}}
\newcommand{\begs}{\begin{examples}}
\newcommand{\eegs}{\end{examples}}
\newcommand{\bdefe}{\begin{defn}}
\newcommand{\edefe}{\end{defn}}

\newcommand{\bca}{\begin{case}}
\newcommand{\eca}{\end{case}}

\newcommand{\bprob}{\begin{prob}}
\newcommand{\eprob}{\end{prob}}
\newcommand{\bques}{\begin{ques}}
\newcommand{\eques}{\end{ques}}
\newcommand{\bei}{\begin{itemize}}
\newcommand{\eei}{\end{itemize}}
\newcommand{\bcl}{\begin{claim}}
\newcommand{\ecl}{\end{claim}}

\newcommand{\bsol}{\begin{mysolution}}
\newcommand{\esol}{\end{mysolution}}

\newcommand{\bcon}{\begin{conj}}
\newcommand{\econ}{\end{conj}}
\newcommand{\bcons}{\begin{conjs}}
\newcommand{\econs}{\end{conjs}}
\newcommand{\bprop}{\begin{prop}}
\newcommand{\eprop}{\end{prop}}
\newcommand{\br}{\begin{rem}}
\newcommand{\er}{\end{rem}}
\newcommand{\brs}{\begin{rems}}
\newcommand{\ers}{\end{rems}}
\newcommand{\bo}{\begin{obser}}
\newcommand{\eo}{\end{obser}}
\newcommand{\bos}{\begin{obsers}}
\newcommand{\eos}{\end{obsers}}
\newcommand{\bpf}{\begin{pf}}
\newcommand{\epf}{\end{pf}}
\newcommand{\ba}{\begin{array}}
\newcommand{\ea}{\end{array}}
\newcommand{\beq}{\begin{eqnarray}}
\newcommand{\beqq}{\begin{eqnarray*}}
\newcommand{\eeq}{\end{eqnarray}}
\newcommand{\eeqq}{\end{eqnarray*}}

\newcommand{\ds}{\displaystyle}

\begin{document}
\bibliographystyle{amsplain}
\title [Schwarz-Pick and Landau type theorems for solutions to the Dirichlet-Neumann problem in the unit disk]
{Schwarz-Pick and Landau type theorems for solutions to the Dirichlet-Neumann problem in the unit disk}
\author{Peijin Li}
\address{ Department of Mathematics,
Hunan First Normal University, Changsha, Hunan 410205, People's Republic of China}
\email{wokeyi99@163.com}

\author{Qinghong Luo}
\address{ Department of Mathematics,
Hunan First Normal University, Changsha, Hunan 410205, People's Republic of China}
\email{luoqh207@qq.com}

\author{Saminathan Ponnusamy
}
\address{S. Ponnusamy, Department of Mathematics,
Indian Institute of Technology Madras, Chennai-600 036, India. }
\email{samy@iitm.ac.in}

\subjclass[2010]{Primary: 30C62, 31A30, 30C80; Secondary:  30C20, 31A05.}

\keywords{Dirichlet-Neumann problem, Schwarz's Lemma, Landau type theorem and univalent function.
}

\begin{abstract}
The aim of this paper is to establish some properties of solutions to the Dirichlet-Neumann problem:
$(\partial_z\partial_{\overline{z}})^2 w=g$ in the unit disc $\ID$, $w=\gamma_0$ and $\partial_{\nu}\partial_z\partial_{\overline{z}}w=\gamma$ on $\mathbb{T}$ (the unit circle),
$\frac{1}{2\pi i}\int_{\mathbb{T}}w_{\zeta\overline{\zeta}}(\zeta)\frac{d\zeta}{\zeta}=c$, where $\partial_\nu$ denotes differentiation in the outward normal direction.
More precisely, we obtain  Schwarz-Pick type inequalities and Landau type theorem for solutions to the Dirichlet-Neumann problem.
\end{abstract}


\maketitle
\pagestyle{myheadings}
\markboth{Peijin Li, Qinghong Luo and Saminathan Ponnusamy
}{Schwarz-Pick and Landau type theorems for the solutions to Dirichlet-Neumann problem}

\section{Introduction}\label{sec-1}


Let $\mathbb{C}\cong\mathbb{R}^2$ denote the complex plane and for $r>0$, let ${\mathbb D}_r=\{z\in \IC:\,|z|<r\}$. Denote by  $\mathbb {D}:= {\mathbb D} _1$,
the open unit disk, $\mathbb{T}=\partial\mathbb{D}$, the boundary of $\mathbb{D}$, and
$\overline{\ID}=\ID\cup \mathbb{T}$, the closure of $\mathbb{D}$. Denote by $\mathcal{C}(\Omega)$, the set of all continuous
functions in a domain $\Omega$ in $\IC$.   The space of integrable functions in $\Omega$ is denoted by $L^{1}(\Omega)$.
Denote by ${\mathcal H}(\ID,\ID)$ (resp. ${\mathcal A}(\ID,\ID)$) the class of all complex-valued harmonic (resp. analytic) self-mappings of the unit disk $\ID$

The Dirichlet and the Neumann boundary value problems in complex analysis have been very well studied in the literature. See
\cite{Be2, BV} for investigations of basic boundary value problems with different kinds of boundary conditions.

In this paper we investigate some properties of solutions to the following Dirichlet-Neumann problem:
\be\label{eq-DN1}
\begin{cases} (\partial_z\partial_{\overline{z}})^2 w=g &
\mbox{ in }\ds  \ID,\\
\ds w=\gamma_0 & \mbox{ on }\ds \mathbb{T},\\
\ds \partial_{\nu}\partial_z\partial_{\overline{z}}w=\gamma & \mbox{ on }\ds \mathbb{T}
\end{cases}
\ee
and
\be\label{eq-DN2}
\frac{1}{2\pi i}\int_{\mathbb{T}}w_{\zeta\overline{\zeta}}(\zeta)\frac{d\zeta}{\zeta}=c,
\ee
where $\partial_\nu$ denotes differentiation in the outward normal direction, $g\in L^1(\ID)$, $\gamma_0$,  $\gamma\in C(\mathbb{T})$,  $c\in\mathbb{C}$ is a constant, and satisfying the condition
\be\label{eq-DN3}
\frac{1}{2\pi}\int^{2\pi}_{0}\gamma(e^{it})\,dt=\frac{2}{\pi}\int_{\ID}g(\zeta)\,d A(\zeta).
\ee

Here
$$\partial _z:=\frac{\partial}{\partial z}=\frac{1}{2} \left ( \frac{\partial}{\partial x}-i\frac{\partial}{\partial y}\right )~\mbox{ and }~
\partial _{\overline{z}}:= \frac{\partial}{\partial \overline{z}}=\frac{1}{2} \left ( \frac{\partial}{\partial x}+i\frac{\partial}{\partial y}\right )
$$
represent the complex differential operators so that $\Delta =4 \partial_z\partial_{\overline{z}}$, is the Laplacian.
So, $\partial_z\partial_{\overline{z}}$ and $(\partial_z\partial_{\overline{z}})^2$ are called the harmonic and biharmonic
operators, respectively. Consequently, $\partial_z\partial_{\overline{z}}w=0$
in $\ID$ is equivalent to the statement that $w$ is harmonic, while a solution to the equation
$(\partial_z\partial_{\overline{z}})^2 w=0$ is called a {\it biharmonic function}.
See \cite{CPW0, CPW4} and the references therein for certain properties of biharmonic functions.

For $z\in\ID$, let
$$P(z, e^{it})=\frac{1-|z|^2}{|1-ze^{-it}|^2}
$$
denote the {\it Poisson kernel} in $\ID$. The function $z\mapsto P(z, e^{it})$ is harmonic in $\ID$.

In \cite{BV}, it was shown that the condition \eqref{eq-DN3} ensures that
all solutions to \eqref{eq-DN1} satisfying the condition \eqref{eq-DN2} are given by the formula
\be\label{eq-DNS}
w(z)=-c(1-|z|^2)+\mathcal{P}_{\gamma_0}(z)+\mathcal{G}_1[\gamma](z)-\mathcal{G}_2[g](z),
\ee
where
$$\mathcal{P}_{\gamma_0}(z)=\frac{1}{2\pi}\int^{2\pi}_{0}P(z, e^{it})\gamma_0(e^{it})\,dt,$$
\be\label{eq-DNS1}
\mathcal{G}_1[\gamma](z)=\frac{1}{4\pi}\int^{2\pi}_{0}H_2(z, e^{it})\gamma(e^{it})\,dt,
\ee
and
\be\label{eq-DNS2}
\mathcal{G}_2[g](z)=\int_{\ID}H_2(z, \zeta)g(\zeta)\,d A(\zeta),
\ee
with
\beq\label{eq-DNS3}
H_2(z, \zeta)&=&-|\zeta-z|^2\log|\zeta-z|^2\nonumber\\
&&-(1-|z|^2)\left [4+\frac{1-z\overline{\zeta}}{z\overline{\zeta}}\log(1-z\overline{\zeta})
+\frac{1-\overline{z}\zeta}{\overline{z}\zeta}\log(1-\overline{z}\zeta)\right]\nonumber\\
&&-\frac{(\zeta-z)(1-z\overline{\zeta})}{z}\log(1-z\overline{\zeta})
-\frac{\overline{(\zeta-z)}(1-\overline{z}\zeta)}{\overline{z}}\log(1-\overline{z}\zeta),
\eeq
and $dA(\zeta) =(1/\pi)\, dx\,dy$ denotes the normalized area measure in $\ID$.

\section{Main Results}\label{sec-2}

\subsection{A Schwarz type lemma}

The classical Schwarz lemma states that if $f\in {\mathcal A}(\ID,\ID)$ such that $f(0)=0$, then
$|f(z)|\leq|z|$ for all $z\in\mathbb{D}$. This result has been a crucial result in many branches of
research for more than a hundred years. The harmonic analog of this statement states that  if $f\in {\mathcal H}(\ID,\ID)$ such that $f(0)=0$,
then (Heniz \cite{He})
\be\label{eq-h1}
|f(z)|\leq\frac{4}{\pi}\arctan|z| ~\mbox{ for $z\in\mathbb{D}$}.
\ee
Later, by removing the assumption $f(0)=0$, Pavlovi\'{c} \cite[Theorem 3.6.1]{P} established the following sharp inequality for $f\in {\mathcal H}(\ID,\ID)$:
\be\label{eq-h2}
\left|f(z)-\frac{1-|z|^2}{1+|z|^2}f(0)\right|\leq\frac{4}{\pi}\arctan|z| ~\mbox{ for $z\in\mathbb{D}$}.
\ee

The first purpose of this paper is to consider results of the above
type for solutions to \eqref{eq-DN1} satisfying the conditions \eqref{eq-DN2} and \eqref{eq-DN3}.

\begin{thm}\label{thm-1}
Suppose that $g\in\mathcal{C}(\overline{\mathbb{D}})$ and $\gamma\in \mathcal{C}(\mathbb{T})$, and that $w\in\mathcal{C}^{4}(\mathbb{D})\bigcap\mathcal{C}(\overline{\mathbb{D}})$
satisfying the equation \eqref{eq-DN1} with the conditions \eqref{eq-DN2} and \eqref{eq-DN3}. Then
for $z\in\overline{\mathbb{D}}$,
\be\label{eqs}
\left|w(z)-\frac{1-|z|^2}{1+|z|^2}\mathcal{P}_{\gamma_0}(0)\right|\leq \frac{4}{\pi}\|\mathcal{P}_{\gamma_0}\|_{\infty}\arctan|z|+|c|
+\|\gamma\|_{\infty}N_1(|z|)+\|g\|_{\infty}N_2(|z|),
\ee
where
$$N_1(t)=2\log4+\frac{1-t^2}{2}\left (\frac{2}{3}\pi^2-2\right )^{\frac{1}{2}}+4\left (\frac{\pi^2}{6}\right )^{\frac{1}{2}},
$$
$$ N_2(t)=4\log4+(1-t^2)\left (\frac{2}{3}\pi^2-2\right )^{\frac{1}{2}}+\frac{16}{3}\left (\frac{\pi^2}{6}\right )^{\frac{1}{2}},
$$
$\|\mathcal{P}_{\gamma_0}\|_{\infty}=\sup_{z\in\mathbb{D}}\{|\mathcal{P}_{\gamma_0}(z)|\}$, $\|\gamma\|_{\infty}=\sup_{z\in\mathbb{T}}\{|\gamma(z)|\}$
and $\|g\|_{\infty}=\sup_{z\in\mathbb{D}}\{|g(z)|\}$.
\end{thm}

Clearly, if $c=0$, $\gamma=g\equiv 0$ and $w$ maps $\ID$ into itself, then \eqref{eqs} coincides with \eqref{eq-h2}.

\subsection{A Schwarz-Pick type lemma}

For a $2\times 2$ real matrix $M:=M_{2\times 2}$,
the matrix norm and the matrix function are defined by
$$\|M\|=\sup\{|Mz|:\,z\in \mathbb{T}\}
~\mbox{ and }~  \lambda(M)=\inf\{|Mz|:\;z\in \mathbb{T}\},
$$
respectively. For a complex-valued function $w = f(z)=u(z)+iv(z)$, the
Jacobian matrix $D_f$ and Jacobian (determinant) $J_f$ of $f$ are defined by
$$D_f = \left(
  \begin{array}{cc}
    u_x & u_y \\
    v_x & v_y \\
  \end{array}
\right) ~\mbox{ and }~ J_f :=\det D_{f} =u_xv_y-v_xu_y =|f_{z}|^{2}-|f_{\overline{z}}|^{2},
$$
respectively. Obviously
\be\label{pj-11}
\|D_f(z)\|=\sup\{|D_f(z)\varsigma|:\,|\varsigma|=1\}=|f_z(z)|+|f_{\overline z}(z)|
\ee
and
$$
\lambda(D_f(z))=\inf\{|D_f(z)\varsigma|:\,|\varsigma|=1\}=\big ||f_z(z)|-|f_{\overline z}(z)|\big |.
$$

Colonna \cite{Co} obtained a sharp Schwarz-Pick type lemma for $f\in {\mathcal H}(\ID,\ID)$:
\be\label{eq-Co}
\|D_{f}(z)\|\leq\frac{4}{\pi}\frac{1}{1-|z|^{2}} ~\mbox{  for $z\in\mathbb{D}$}.
\ee

Our second aim in this paper is to prove the following Schwarz-Pick type lemma for solutions to \eqref{eq-DN1} satisfying the conditions
\eqref{eq-DN2} and \eqref{eq-DN3}.

\begin{thm}\label{thm-2}
Suppose that $g\in\mathcal{C}(\overline{\mathbb{D}})$ and $\gamma\in \mathcal{C}(\mathbb{T})$, and that $w\in\mathcal{C}^{4}(\mathbb{D})\bigcap\mathcal{C}(\overline{\mathbb{D}})$
satisfying the equation \eqref{eq-DN1} with the conditions \eqref{eq-DN2} and \eqref{eq-DN3}. Then
for $z\in\mathbb{D}$,
\be\label{eqsp}
\|D_{w}(z)\|\leq \frac{4}{\pi}\|\mathcal{P}_{\gamma_0}\|_{\infty}\frac{1}{1-|z|^2}+2|c|+\|\gamma\|_{\infty}N_3(|z|)+\|g\|_{\infty}N_4(|z|),
\ee
where
$$N_3(t)=2\left (\frac{\pi^2}{3}+1\right )^{\frac{1}{2}}+t\left(\frac{2}{3} \pi^2 -2\right)^\frac{1}{2}
+(1-t^2)\left (\frac{\pi^2}{6}-1\right )^{\frac{1}{2}} + \left(\frac{\pi^2}{3}-\frac{1}{2}\right )^{\frac{1}{2}},$$
$$N_4(t)=2(\log4+1)+t\left ( \frac{2}{3}\pi^2-2\right )^{\frac{1}{2}}
+\frac{2}{3}(1-t^2)\left (\frac{\pi^2}{6}-1\right)^{\frac{1}{2}} +\frac{2}{3}\left(\frac{\pi^2}{3}-\frac{1}{2}\right)^{\frac{1}{2}},$$
 $\|\mathcal{P}_{\gamma_0}\|_{\infty}$,  $\|\gamma\|_{\infty}$ and  $\|g\|_{\infty}$ are as in Theorem {\rm \ref{thm-1}}.
\end{thm}

\br
We note that if $c=0$, $\gamma=g \equiv 0$ and $w$ maps $\ID$ into itself, then \eqref{eqsp} coincides with \eqref{eq-Co}.
\er
\subsection{A Landau type theorem}

The classical Landau theorem says that there is a $\rho=\frac{1}{M+\sqrt{M^{2}-1}}$ such that every function $f$, analytic in $\mathbb{D}$ with $f(0)=f^{'}(0)-1=0$ and $|f(z)|< M$, is univalent in the disk $\mathbb{D}_{\rho}$. Moreover, the range $f(\mathbb{D}_{\rho})$ contains a disk of radius $M\rho^{2}$, where $M\geq1$ is a constant (see \cite{E}).
The Landau theorem has become an important tool in geometric function theory of
one complex variable (cf. \cite{Br, Za}). Unfortunately, for a general class of functions, there is no Landau type theorem (cf. \cite{GP, Wu}).
To establish analogs of the Landau type theorem for more general classes of functions, it is necessary to restrict our focus to certain subclasses
(cf. \cite{AM, CGH, CG, CLW, CPW0, CPW5, CV, LWX, Wu}).

As an application of Theorem \ref{thm-2}, we get the following Landau type theorem for the solutions to \eqref{eq-DN1} with the conditions \eqref{eq-DN2} and \eqref{eq-DN3}.

\begin{thm}\label{thm-3}
Suppose that $\gamma\in C(\mathbb{T})$, $g\in C(\overline{\ID})$, that $w\in\mathcal{C}^{4}(\mathbb{D})\bigcap\mathcal{C}(\overline{\mathbb{D}})$
satisfying the equation \eqref{eq-DN1} with the conditions \eqref{eq-DN2}, \eqref{eq-DN3} and $w(0)=J_w(0)-1=0$, and that
$\|\gamma_0\|_{\infty}\leq L_1$, $\|\gamma\|_{\infty}\leq L_2$ and $\|g\|_{\infty}\leq L_3$, where $L_j, j\in\{1, 2, 3\}$, are constants.
Then
\begin{enumerate}
  \item \label{eq1.7}
  $w$ is univalent in $\mathbb{D}_{r_{0}}$, where $r_{0}$ satisfies the following equation
\be\label{fil-1}
\frac{1}{L_4}-2r_0\Big(\frac{4L_1}{\pi}\frac{2-r_0}{(1-r_0)^2}+L_5\Big)-8L_3\log\frac{1+r_0}{1-r_0}=0;
\ee
\item \label{eq1.10}
  $w(\mathbb{D}_{r_{0}})$ contains a univalent disk $\mathbb{D}_{R_{0}}$ with
$$
R_{0}\geq\frac{8L_1}{\pi}\,\Big(\frac{r_{0}}{1-r_{0}}\Big)^2+L_5r_0^2+\frac{8L_3r_0^2(3-r_0^2)}{(1-r_0^2)^2},
$$
where
$$L_4=2|c|+\frac{4}{\pi}L_1+L_2N_3(0)+L_3N_4(0),\,\,\,\,L_5=|c|+L_2M_1+L_3M_2,
$$
$N_3(0)$ and $N_4(0)$ are defined in Theorem \ref{thm-2}, whereas $M_1$ and $M_2$ are defined in Lemmas \ref{lem-1} and \ref{lem-2}, respectively.
\end{enumerate}
\end{thm}

We would like to remark that this article continues the earlier study on this topic from \cite{LP}. Recently, many authors have studied the
Schwarz type lemma,  Schwarz-Pick type lemma and Landau type theorem for solutions of different equations (cf. \cite{CK, CLW, CZ, LRW, LWX}).

The rest of this article is organized as follows. Section \ref{sec-3} is devoted to stating and proving
several useful lemmas. In Section \ref{sec-4}, we present the proofs of Theorems \ref{thm-1}, \ref{thm-2} and \ref{thm-3}.

\section{Several Basic Lemmas}\label{sec-3}

In this section, we shall prove three lemmas which will be used later on.

\begin{lem}\label{lem-1}
Suppose $\gamma\in C(\mathbb{T})$, and $\mathcal{G}_1[\gamma]$ is defined by \eqref{eq-DNS1}. Then for $z\in \ID$,
$$\left|\frac{\partial}{\partial z}\mathcal{G}_1[\gamma](z)-\frac{\partial}{\partial z}\mathcal{G}_1[\gamma](0)\right|\leq
(\|\gamma\|_{\infty} M_1)|z|
$$
and
$$\left|\frac{\partial}{\partial \overline{z}}\mathcal{G}_1[\gamma](z)-\frac{\partial}{\partial \overline{z}}\mathcal{G}_1[\gamma](z)\right|
\leq (\|\gamma\|_{\infty} M_1)|z|,
$$
where
\beqq
M_1=\frac{1}{2}\left[\frac{2\pi}{\sqrt{3}}+1+\left(\frac{2}{3} \pi^2 -2\right)^\frac{1}{2}
+\left (\frac{\pi^2}{6}-\frac{5}{4}\right )^{\frac{1}{2}}
+\left(\frac{\pi^2}{6}-1\right )^{\frac{1}{2}}
+\left(\frac{\pi^2}{3}-\frac{11}{4}\right )^{\frac{1}{2}}\right].
\eeqq
\end{lem}

\bpf
To prove the lemma, we only need to show the first inequality, namely,
$$\left|\frac{\partial}{\partial z}\mathcal{G}_1[\gamma](z)-\frac{\partial}{\partial z}\mathcal{G}_1[\gamma](0)\right|\leq (\|\gamma\|_{\infty} M_1)|z|,
$$
since the proof of the second inequality is similar. Let
\beqq
I_1(z)&=&\frac{1}{2\pi}\int^{2\pi}_{0}(e^{-it}-\overline{z})\big(\log|1-ze^{-it}|^2+1\big)\gamma(e^{it})\,dt,\\
I_2(z)&=&\frac{1}{2\pi}\int^{2\pi}_{0}\overline{z}\left[4+\frac{1-ze^{-it}}{ze^{-it}}\log(1-ze^{-it})
+\frac{1-\overline{z}e^{it}}{\overline{z}e^{it}}\log(1-\overline{z}e^{it})\right]\gamma(e^{it})\,dt,\\
I_3(z)&=& \frac{1}{2\pi}\int^{2\pi}_{0}(1-|z|^2)\left(\frac{\log(1-ze^{-it})}{z^2e^{-it}}+\frac{1}{z}\right)\gamma(e^{it})\,dt,~\mbox{ and }~\\
I_4(z)&=& \frac{1}{2\pi}\int^{2\pi}_{0}\left(\frac{e^{it}-z^2e^{-it}}{z^2}\log(1-ze^{-it})+\frac{1-ze^{-it}}{z}\right)\gamma(e^{it})\,dt.
\eeqq
Now, we need to estimate $|I_j(z)-I_j(0)|$ for $j=1, 2, 3, 4$, respectively. 
\bcl \label{eq-I1}
$\ds |I_1(z)-I_1(0)|\leq \|\gamma\|_{\infty}|z|\left[2\left (\frac{\pi^2}{3}\right )^{\frac{1}{2}}+1\right].$
\ecl

By elementary calculations, we get
$$|I_1(z)-I_1(0)|\leq \|\gamma\|_{\infty}(2J_1(z)+|z|),$$
where $$J_1(z)=\frac{1}{2\pi}\int^{2\pi}_{0}\left|\log|1-ze^{-it}|^2\right| \,dt
=\frac{1}{2\pi}\int^{2\pi}_{0}\left|\sum^{\infty}_{n=1}\frac{z^ne^{-int}}{n}+\sum^{\infty}_{n=1}\frac{\overline{z}^ne^{int}}{n}\right| \,dt.$$
By the H\"{o}lder inequality and Parseval's theorem, we obtain
\beqq
J_1(z)
&\leq & \left(\frac{1}{2\pi}\int^{2\pi}_{0}\left|\sum^{\infty}_{n=1}\frac{z^ne^{-int}}{n}+\sum^{\infty}_{n=1}\frac{\overline{z}^ne^{int}}{n}\right|^2 dt\right)^{\frac{1}{2}}\\
&=& \left(2\sum^{\infty}_{n=1}\frac{|z|^{2n}}{n^2}\right)^{\frac{1}{2}}
\leq |z|\left (\frac{\pi^2}{3}\right )^{\frac{1}{2}}
\eeqq
which proves Claim \ref{eq-I1}.


\bcl\label{eq-I2}
$\ds  |I_2(z)-I_2(0)|\leq\|\gamma\|_{\infty}|z|\left(\frac{2}{3} \pi^2 -2\right)^\frac{1}{2}.$
\ecl

Since $|I_2(z)-I_2(0)|=|I_2(z)|$, Claim \ref{eq-I2} follows from \cite[Claim 2.4]{LP}.


\bcl\label{eq-I3}
$\ds  |I_3(z)-I_3(0)|\leq \|\gamma\|_{\infty}|z|\left[\left (\frac{\pi^2}{6}-\frac{5}{4}\right )^{\frac{1}{2}}+\left(\frac{\pi^2}{6}-1\right )^{\frac{1}{2}}\right].$
\ecl

Since
$$\frac{\log (1-ze^{-it})}{z^2e^{-it}} +\frac{1}{z}= -\sum^{\infty}_{n=1}\frac{z^{n-1}e^{-int}}{n+1},$$
we deduce that
$$I_3(0)=\frac{1}{2\pi}\int^{2\pi}_{0}\left(-\frac{e^{-it}}{2}\right)\gamma(e^{it})\,dt.
$$
Then
$$|I_3(z)-I_3(0)|\leq \|\gamma\|_{\infty}\big(J_2(z)+J_3(z)\big),$$
where
$$J_2(z)=\frac{1}{2\pi}\int^{2\pi}_{0}\left|\sum^{\infty}_{n=2}\frac{z^{n-1}e^{-int}}{n+1}\right|\,dt
~\mbox{and}~
J_3(z)=\frac{|z|^2}{2\pi}\int^{2\pi}_{0}\left|\sum^{\infty}_{n=1}\frac{z^{n-1}e^{-int}}{n+1}\right|\,dt.$$
As in the proof of the estimate for $J_1(z)$, by the H\"{o}lder inequality and Parseval's theorem, we obtain
\beqq
J_2(z)
&\leq & \left(\frac{1}{2\pi}\int^{2\pi}_{0}\left|\sum^{\infty}_{n=2}\frac{z^{n-1}e^{-int}}{n+1}\right|^2 dt\right)^{\frac{1}{2}}\\
&=& \left(\sum^{\infty}_{n=2}\frac{|z|^{2(n-1)}}{(n+1)^2}\right)^{\frac{1}{2}}
\leq |z|\left (\frac{\pi^2}{6}-\frac{5}{4}\right )^{\frac{1}{2}}.
\eeqq
Similarly, we know that
$$J_3(z)\leq |z|^2\left (\frac{\pi^2}{6}-1\right )^{\frac{1}{2}}.
$$
Claim \ref{eq-I3} follows as $|z|<1$.

\bcl\label{eq-I4}
$\ds  |I_4(z)-I_4(0)|\leq\|\gamma\|_{\infty}|z|\left(\frac{\pi^2}{3}-\frac{11}{4}\right )^{\frac{1}{2}}.$
\ecl

We rewrite $I_4(z)$ as
$$I_4(z) =\frac{1}{2\pi}\int^{2\pi}_{0} e^{-it}G(z e^{-it}) \gamma(e^{it})\,dt,
$$
where
$$G(z)=\frac{1-z^2}{z^2}\log (1-z) +\frac{1}{z}-1= 2\sum^{\infty}_{n=1}\frac{z^{n}}{n(n+2)}-\frac{3}{2}.
$$
Then $I_4(0) =\frac{1}{2\pi}\int^{2\pi}_{0} (-\frac{3}{2})e^{-it} \gamma(e^{it})\,dt$.
As before, we find that
\beqq
|I_4(z)-I_4(0)|&\leq&\frac{\|\gamma\|_{\infty}}{2\pi}\int^{2\pi}_{0}\left|2\sum^{\infty}_{n=1}\frac{z^{n}e^{-int}}{n(n+2)}\right| dt\\
&\leq& \|\gamma\|_{\infty}\left(\frac{1}{2\pi}\int^{2\pi}_{0}\left| 2\sum^{\infty}_{n=1}\frac{z^{n}e^{-int}}{n(n+2)}\right|^2 dt\right)^{\frac{1}{2}}\\
& = &\|\gamma\|_{\infty} \left( \sum^{\infty}_{n=1}\frac{4|z|^{2n}}{n^2(n+2)^2} \right)^{\frac{1}{2}}
\leq\|\gamma\|_{\infty}|z|\left(\frac{\pi^2}{3}-\frac{11}{4}\right )^{\frac{1}{2}},
\eeqq
since
$$ \sum^{\infty}_{n=1}\frac{4}{n^2(n+2)^2}=\sum^{\infty}_{n=1}\frac{1}{(n+2)^2} +
\sum^{\infty}_{n=1}\frac{1}{n^2}-2\sum^{\infty}_{n=1}\frac{1}{n(n+2)} =
\frac{\pi^2}{3}-\frac{11}{4},
$$
and Claim \ref{eq-I4} follows.

Therefore, by Claims \ref{eq-I1}, \ref{eq-I2}, \ref{eq-I3}, \ref{eq-I4} and \cite[Proposition 2.4]{KP},
we conclude that
$$
\left|\frac{\partial}{\partial z}\mathcal{G}_1[\gamma](z)-\frac{\partial}{\partial z}\mathcal{G}_1[\gamma](0)\right|
\leq\frac{1}{2}\sum^{4}_{j=1}|I_j(z)-I_j(0)|\leq (\|\gamma\|_{\infty} M_1)|z|,
$$
as required.
\epf

\begin{lem}\label{lem-2}
Suppose $g\in C(\overline{\ID})$ and $\mathcal{G}_2[g]$ is defined in \eqref{eq-DNS2}. Then for $z\in \ID_{r^{\star}}$,
$$\left|\frac{\partial}{\partial z}\mathcal{G}_2[g](z)-\frac{\partial}{\partial z}\mathcal{G}_2[g](0)\right|
\leq \|g\|_{\infty}|z|M_2+4\|g\|_{\infty}\log\frac{1+|z|}{1-|z|}$$
and
$$\left|\frac{\partial}{\partial \overline{z}}\mathcal{G}_2[g](z)-\frac{\partial}{\partial \overline{z}}\mathcal{G}_2[g](0)\right|\leq \|g\|_{\infty}|z|M_2+4\|g\|_{\infty}\log\frac{1+|z|}{1-|z|},$$
where $0\leq r^{\star}<2r_0$ and $r_0$ is determined by Eqn. \eqref{fil-1}, and
\beqq
M_2=\log4+1
+\left ( \frac{2}{3}\pi^2-2\right )^{\frac{1}{2}}
+\frac{2}{3}\left[\left ( \frac{\pi^2}{6} -\frac{5}{4}\right )^{\frac{1}{2}}+\left (\frac{\pi^2}{6}-1\right)^{\frac{1}{2}} +\left(\frac{\pi^2}{3}-\frac{11}{4}\right)^{\frac{1}{2}}\right].
\eeqq
\end{lem}

\bpf
To prove the two inequalities, we only need to prove the first inequality, namely,
$$\left|\frac{\partial}{\partial z}\mathcal{G}_2[g](z)-\frac{\partial}{\partial z}\mathcal{G}_2[g](0)\right|\leq
(\|g\|_{\infty} M_2)|z| ,
$$
because the proof of the second inequality is similar. To do this, we let
\beqq
I_5(z)&=&\int_{\ID}(\overline{\zeta}-\overline{z})\big(\log|\zeta-z|^2+1\big)g(\zeta)\,d A(\zeta),\\
I_6(z)&=&\int_{\ID}\overline{z}\left[4+\frac{1-z\overline{\zeta}}{z\overline{\zeta}}\log(1-z\overline{\zeta})
+\frac{1-\overline{z}\zeta}{\overline{z}\zeta}\log(1-\overline{z}\zeta)\right]g(\zeta)\,d A(\zeta),\\
I_7(z)&=&\int_{\ID}(1-|z|^2)\left(\frac{\log(1-z\overline{\zeta})}{z^2\overline{\zeta}}+\frac{1}{z}\right)g(\zeta)\,d A(\zeta) ,~\mbox{ and }~\\
I_8(z)&=&\int_{\ID}\left(\frac{\zeta-z^2\overline{\zeta}}{z^2}\log(1-z\overline{\zeta})+\frac{|\zeta|^2-z\overline{\zeta}}{z}\right)g(\zeta)\,d A(\zeta).
\eeqq
In the following, we estimate $|I_j(z)-I_j(0)|$ for $j=5, 6, 7, 8$, respectively.
\bcl\label{eq-I5}
$\ds  |I_5(z)-I_5(0)|\leq \|g\|_{\infty}\Big[4\log\frac{1+|z|}{1-|z|}+(\log4+1)|z|\Big].$
\ecl

By calculations, we get
$$I_5(z)-I_5(0)=-J_4(z)+J_5(z),$$
where
$$J_4(z)=\int_{\ID}\overline{z}\big(\log|\zeta-z|^2+1\big)g(\zeta)\,d A(\zeta)$$
and
$$J_5(z)=\int_{\ID}\overline{\zeta}\big(\log|\zeta-z|^2-\log|\zeta|^2\big)g(\zeta)\,d A(\zeta).$$
Obviously,
\beq\label{eq-J4}|J_4(z)|\leq (\log4+1)\|g\|_{\infty}|z|. \eeq
In order to estimate $J_5(z)$, we let
$$h(z, \zeta)=\overline{\zeta}\log|\zeta-z|^2.$$
Then, for $z\in \ID_{r^{\star}}$, by Fubini's Theorem, we get
\beqq
|J_5(z)|&\leq& \|g\|_{\infty}\int_{\ID}\left|h(z, \zeta)-h(0, \zeta)\right|\,d A(\zeta)\\
&\leq&\|g\|_{\infty}\int_{\ID}\left(\int_{[0,z]}\Big|h_z(z, \zeta)dz+h_{\overline{z}}(z, \zeta)d\overline{z}\Big|\right)\,d A(\zeta)\\
&\leq&\|g\|_{\infty}\int_{\ID}\left(\int_{[0,z]}\Big(|h_z(z, \zeta)|+|h_{\overline{z}}(z, \zeta)|\Big)|dz|\right)\,d A(\zeta)\\
&=&\|g\|_{\infty}\int_{[0,z]}H(z, \zeta)|dz|,
\eeqq
where
$$H(z, \zeta)=\int_{\ID}\Big(|h_z(z, \zeta)|+|h_{\overline{z}}(z, \zeta)|\Big)\,d A(\zeta)=2\int_{\ID}\left|\frac{\zeta}{\zeta-z}\right|\,d A(\zeta).$$
In order to estimate $H(z, \zeta)$, we let
$$\zeta \mapsto \eta =\phi (\zeta)=\frac{z-\zeta}{1-\zeta\overline{z}}=re^{i\theta}
$$
so that $\phi =\phi^{-1}$,
$$\zeta=\frac{z-\eta}{1-\eta\overline{z}},~~z-\zeta=\frac{\eta(1-|z|^2)}{1-\eta\overline{z}}, ~~\phi ' (\zeta)=-\frac{1-|z|^2}{(1-\zeta\overline{z})^2},
$$
and thus,
$$d A(\zeta)= |({\phi^{-1}}) ' (\eta)|^2 d A(\eta)=\frac{(1-|z|^2)^2}{|1-\eta\overline{z}|^4}d A(\eta).
$$
Consequently, switching to polar coordinates yields
\beqq
H(z, \zeta)&=&\int_{\ID}\frac{|z-\eta|(1-|z|^2)}{|\eta|\cdot|1-\eta\overline{z}|^4}\,d A(\eta)
\leq\frac{2(1-|z|^2)}{\pi}\int^{1}_{0}\int^{2\pi}_{0}\frac{1}{|1-\overline{z}re^{i\theta}|^4}\,d\theta\, dr.
\eeqq
By Parseval's theorem, we get
$$\frac{1}{2\pi}\int^{2\pi}_{0}\frac{d\theta}{|1-\overline{z}re^{-i\theta}|^{4}}=\sum^{\infty}_{n=0}(n+1)^2|z|^{2n}r^{2n},$$
and thus,
$$H(z, \zeta)\leq 4(1-|z|^2)\sum^{\infty}_{n=0}\frac{(n+1)^2}{2n+1}|z|^{2n}\leq 4(1-|z|^2)\sum^{\infty}_{n=0}(n+1)|z|^{2n}=\frac{4}{1-|z|^2},$$
because $\sum^{\infty}_{n=0}(n+1)z^{n} =1/(1-z)^2$ for all $|z|<1$. Hence,
\beq\label{eq-J5}
|J_5(z)|\leq \|g\|_{\infty}\int_{[0,z]}\frac{8}{1-|z|^2}|dz|=4\|g\|_{\infty}\log\frac{1+|z|}{1-|z|}.
\eeq
Claim \ref{eq-I5} follows from $|I_5(z)-I_5(0)|\leq|J_4(z)|+|J_5(z)|$, \eqref{eq-J4} and \eqref{eq-J5}.

In the following, we let $\zeta=\rho e^{it}$. By calculations and the H\"{o}lder inequality, we get

\bcl\label{eq-I6}
$\ds  |I_6(z)-I_6(0)|\leq \|g\|_{\infty}|z|\left ( \frac{2}{3}\pi^2-2\right )^{\frac{1}{2}}.$
\ecl

Claim \ref{eq-I6} follows from \cite[Claim 2.8]{LP}, because $|I_6(z)-I_6(0)|=|I_6(z)|$.

\bcl\label{eq-I7}
$\ds  |I_7(z)-I_7(0)|\leq \frac{2}{3}\|g\|_{\infty}|z|\left[\left ( \frac{\pi^2}{6} -\frac{5}{4}\right )^{\frac{1}{2}}+\left ( \frac{\pi^2}{6} -1\right )^{\frac{1}{2}}\right].$
\ecl

As in the proof of Claim \ref{eq-I3} in Lemma \ref{lem-1}, we use the representation
$$\frac{\log(1-z\overline{\zeta})}{z^2\overline{\zeta}}+\frac{1}{z}=-\overline{\zeta}\sum^{\infty}_{n=1}\frac{(z \overline{\zeta})^{n-1}}{n+1}
$$
and obtain that
$$|I_7(z)-I_7(0)|\leq \|g\|_{\infty}\big(J_6(z)+J_7(z)\big),$$
where
$$J_6(z)=\int_{\ID}\left |\sum^{\infty}_{n=2}\frac{(z \overline{\zeta})^{n-1}}{n+1}\right | |\zeta|\,d A(\zeta)
~\mbox{and}~
J_7(z)=|z|^2\int_{\ID}\left |\sum^{\infty}_{n=1}\frac{(z \overline{\zeta})^{n-1}}{n+1}\right | |\zeta|\,d A(\zeta).$$
By H\"{o}lder's inequality and Parseval's theorem, we get
\beqq
J_6(z)&\leq& 2\int^{1}_{0}\left(\frac{1}{2\pi}\int^{2\pi}_{0}\left|\sum^{\infty}_{n=2}\frac{(z\rho e^{-it})^{n-1}}{n+1}\right|^2 dt\right)^{\frac{1}{2}}\rho^2\,d\rho\\
&=& 2\int^{1}_{0}\left(\sum^{\infty}_{n=1}\frac{|z|^{2n}\rho^{2n}}{(n+2)^2}\right)^{\frac{1}{2}}\rho^2\,d\rho\\
&\leq&2|z|\left(\sum^{\infty}_{n=1}\frac{1}{(n+2)^2}\right)^{\frac{1}{2}}   \int^{1}_{0} \rho^2\,d\rho
=\frac{2}{3}|z| \left ( \frac{\pi^2}{6} -\frac{5}{4}\right )^{\frac{1}{2}}.
\eeqq
By similar reasoning as above, one obtains that
$$J_7(z)\leq\frac{2}{3}|z| \left ( \frac{\pi^2}{6} -1\right )^{\frac{1}{2}}.$$
The Claim \ref{eq-I7} follows.

\bcl\label{eq-I8}
$\ds |I_8(z)-I_8(0)|\leq \frac{2}{3}\|g\|_{\infty}|z|\left (\frac{\pi^2}{3}-\frac{11}{4}\right )^{\frac{1}{2}}.
$
\ecl

It follows from
$$\frac{\zeta-z^2\overline{\zeta}}{z^2}\log(1-z\overline{\zeta})+\frac{|\zeta|^2-z\overline{\zeta}}{z}
=\overline{\zeta}\left(\sum^{\infty}_{n=1}\frac{n(1-|\zeta|^2)+2}{n(n+2)}(z \overline{\zeta})^{n} -1-\frac{|\zeta|^2}{2}
\right)$$
that
$$|I_8(z)-I_8(0)|=\left|\int_{\ID}\overline{\zeta}\sum^{\infty}_{n=1}\frac{n(1-|\zeta|^2)+2}{n(n+2)}(z \overline{\zeta})^{n}
g(\zeta)\,d A(\zeta)\right|\leq \|g\|_{\infty} J_8(z),$$
where
$$J_8(z)=\int_{\ID}\left|\sum^{\infty}_{n=1}\frac{n(1-|\zeta|^2)+2}{n(n+2)}(z \overline{\zeta})^{n}\right|\cdot|\zeta|
\,d A(\zeta).$$
Now it is easy to see that 
\beqq
J_8(z)
&\leq&2\int_{0}^{1}\Psi(z, \rho)\rho ^2\, d \rho, \quad \zeta =\rho e^{it},
\eeqq
where
\beqq
\Psi(z, \rho)&=&\left(\frac{1}{2 \pi} \int_{0}^{2 \pi}
\left | \sum^{\infty}_{n=1}\frac{n(1-\rho ^2)+2}{n(n+2)}(z \rho e^{-it})^{n} \right |^{2}\, dt\right)^{\frac{1}{2}}\\
&=&\left (\sum^{\infty}_{n=1}\frac{(n(1-\rho ^2)+2)^2  \rho^{2n}}{n^2(n+2)^2}|z|^{2n}\right )^{\frac{1}{2}}\\
&\leq & |z|\left(4\sum^{\infty}_{n=1}\frac{1}{n^2(n+2)^2}\right )^{\frac{1}{2}}
=|z|\left(\frac{\pi^2}{3}-\frac{11}{4}\right )^{\frac{1}{2}}.
\eeqq
Here we have used the fact that for each $n\geq 1$, $\varphi (\rho)=(n(1-\rho ^2)+2)^2  \rho^{2n}$ is an increasing function $\rho$.
The  Claim \ref{eq-I8} follows.


Now, Claims \ref{eq-I5}, \ref{eq-I6}, \ref{eq-I7}, \ref{eq-I8} and \cite[Proposition 2.4]{KP} guarantee that
$$
\left|\frac{\partial}{\partial z}\mathcal{G}_2[g](z)-\frac{\partial}{\partial z}\mathcal{G}_2[g](0)\right|
\leq\sum^{8}_{j=5}|I_j(z)-I_j(0)|\leq \|g\|_{\infty}|z|M_2+4\|g\|_{\infty}\log\frac{1+|z|}{1-|z|},
$$
as required.
\epf

\begin{lem}\label{lem-3}
For constants $C_j>0$, $j\in\{1, 2, 3, 4\}$, let
$$
\varphi(x)=C_1-C_2x\,\left[\frac{2-x}{(1-x)^{2}}+C_3\right]-C_4\log\frac{1+x}{1-x}, \quad x\in [0,1).
$$
Then  we have
\ben
\item
 $\varphi$ is continuous and strictly decreasing in $(0,1)$;
 \item
 there is a unique $x_{0}\in (0,1)$ such that $\varphi(x_{0})=0$.
 \een
\end{lem}
\bpf
For $x\in[0,1)$, we find that
$$
\varphi'(x)=-\frac{2C_2}{(1-x)^{3}}-C_2C_3-\frac{2C_4}{1-x^2}<0
$$
showing that $\varphi(x)$ is strictly decreasing in $[0,1)$. Moreover,
$$
\varphi(0)=C_1>0 \;\;\text{and}\;\; \lim_{x\to1^{-}}\varphi(x)=-\infty<0
$$
which implies that there is a unique $x_{0}\in (0,1)$ such that $\varphi(x_{0})=0$. The proof of the lemma is complete.
\epf

\section{The proof of main results} \label{sec-4}

In this section, we supply the proofs of  Theorems \ref{thm-1}, \ref{thm-2} and \ref{thm-3}.

\subsection*{ Proof of  Theorem \ref{thm-1}}
By \eqref{eq-DNS} and \eqref{eq-h2}, we can quickly deduce that
$$\left|w(z)-\frac{1-|z|^2}{1+|z|^2}\mathcal{P}_{\gamma_0}(0)\right|
\leq \frac{4}{\pi}\|\mathcal{P}_{\gamma_0}\|_{\infty}\arctan|z|+|c|+|\mathcal{G}_1[\gamma](z)|+|\mathcal{G}_2[g](z)|.
$$
We just need  to estimate $|\mathcal{G}_1[\gamma](z)|$ and $|\mathcal{G}_2[g](z)|$.

\bcl\label{eq-g1}
$\ds |\mathcal{G}_1[\gamma](z)|\leq \|\gamma\|_{\infty}N_1(|z|).
$
\ecl

By elementary calculations, \eqref{eq-DNS1} and Claim \ref{eq-I2}, we obtain
$$|\mathcal{G}_1[\gamma](z)|\leq \frac{\|\gamma\|_{\infty}}{2}\left[4\log4+(1-|z|^2)\left (\frac{2}{3}\pi^2-2\right )^{\frac{1}{2}}\right] +4\|\gamma\|_{\infty}J_9(z),
$$
where
$$J_9(z)=\frac{1}{2 \pi} \int_{0}^{2 \pi}\left|\frac{\log(1-ze^{-it})}{z}\right|\, dt.
$$
Now, in order to estimate $J_9(z)$, we use H\"{o}lder's inequality and Parseval's theorem  to get
\beqq
J_9(z)\leq \left (\frac{1}{2 \pi} \int_{0}^{2 \pi}\left|\sum^{\infty}_{n=1}\frac{z^{n-1}e^{-int}}{n}\right|^2\, dt\right )^{\frac{1}{2}}
=\left (\sum^{\infty}_{n=1}\frac{|z|^{2(n-1)}}{n^2}\right )^{\frac{1}{2}}\leq \left (\frac{\pi^2}{6}\right )^{\frac{1}{2}}
\eeqq
which proves Claim \ref{eq-g1}.

\bcl\label{eq-g2}
$\ds |\mathcal{G}_2[g](z)|\leq \|g|_{\infty}N_2(|z|).
$
\ecl

By elementary calculations, \eqref{eq-DNS2} and Claim \ref{eq-I6}, we obtain
$$|\mathcal{G}_2[g](z)|\leq \|g\|_{\infty}\left[4\log4+(1-|z|^2)\left (\frac{2}{3}\pi^2-2\right )^{\frac{1}{2}}\right]
+8\|g\|_{\infty}J_{10}(z),$$
where
$$J_{10}(z)=\int_{\ID}\left|\frac{\log(1-z\overline{\zeta})}{z}\right|\, d A(\zeta).$$
In order to estimate $J_{10}(z)$, we let $\zeta=\rho e^{it}$. Switching to polar coordinates and by H\"{o}lder's
inequality and Parseval's theorem, we get
\beqq
J_{10}(z)&\leq &2\int_{0}^{1}\left(\frac{1}{2 \pi} \int_{0}^{2 \pi}
\left | \sum^{\infty}_{n=1}\frac{(z \rho e^{-it})^{n-1}}{n} \right |^{2}\, dt\right)^{\frac{1}{2}}\rho ^2\, d \rho\\
&\leq& 2\left(\frac{\pi^2}{6}\right)^{\frac{1}{2}} \int_{0}^{1}\rho ^2\, d \rho
=\frac{2}{3}\left(\frac{\pi^2}{6}\right)^{\frac{1}{2}}
\eeqq
which proves Claim \ref{eq-g2}.

Hence, it follows from Claim \ref{eq-g1} and Claim \ref{eq-g2} that \eqref{eqs} holds, and the proof of the theorem is complete.
\hfill $\Box$

\subsection*{Proof of  Theorem \ref{thm-2}}

By \eqref{eq-DNS} and \eqref{pj-11}, for each $z\in\ID$, we get
$$\|D_w(z)\|=|w_z(z)|+|w_{\overline{z}}(z)|\leq 2|c|+\|D_{\mathcal{P}_{\gamma_0}}(z)\|+\|D_{\mathcal{G}_1[\gamma]}(z)\|+\|D_{\mathcal{G}_2[g]}(z)\|.$$
From \eqref{eq-Co} and \cite[Lemmas 2.2 and 2.3]{LP}, we deduce that
$$\|D_{w}(z)\|\leq \frac{4}{\pi}\|\mathcal{P}_{\gamma_0}\|_{\infty}\frac{1}{1-|z|^2}+2|c|+\|\gamma\|_{\infty}N_3(|z|)+\|g\|_{\infty}N_4(|z|)$$
as required, and the proof of the theorem is complete. \hfill $\Box$

\vspace{7pt}

Before we prove Theorem \ref{thm-3}, let us recall the following result.

\begin{Thm}\label{pj-3}\cite[Lemma 1]{CPW5}
Suppose $f$ is a harmonic mapping of $\ID$ into $\mathbb{C}$ such that $|f(z)|\leq M$ and
$f(z)=\sum^{\infty}_{n=0}a_nz^n+\sum^{\infty}_{n=1}\overline{b}_n\overline{z}^n.$
Then $|a_0|\leq M$ and for all $n\geq 1$,
$$|a_n|+|b_n|\leq\frac{4M}{\pi}.
$$
This estimate is sharp, and the extreme function is
$$ f_n(z)=\begin{cases}
\ds \frac{2M\alpha_1}{\pi}\arg\left(\frac{1+\beta_1z^n}{1-\beta_1z^n}\right),\;|\alpha_1|=|\beta_1|=1, & \ds \mbox{ if }~ n\geq 1,\\
\ds M & \ds \mbox{ if }~   n=0.
\end{cases}
$$
\end{Thm}

\subsection*{ Proof of  Theorem \ref{thm-3}}
The function $\mathcal{P}_{\gamma_0}$ is harmonic in $\ID$ and thus, it can be written in the form
$$\mathcal{P}_{\gamma_0}(z)=\sum^{\infty}_{n=0}a_nz^n+\sum^{\infty}_{n=1}\overline{b}_n\overline{z}^n.$$
Applying Theorem \ref{pj-3}, we get

\vspace{8pt}

$\ds  \left|\frac{\partial}{\partial z}\mathcal{P}_{\gamma_0}(z)-\frac{\partial}{\partial z}\mathcal{P}_{\gamma_0}(0)\right|
+\left|\frac{\partial}{\partial \overline{z}}\mathcal{P}_{\gamma_0}(z)-\frac{\partial}{\partial \overline{z}}\mathcal{P}_{\gamma_0}(0)\right| $
\beq\label{eq-P}
 &=&\left|\sum^{\infty}_{n=2}na_nz^{n-1}\right|+\left|\sum^{\infty}_{n=2}nb_n\overline{z}^{n-1}\right|\nonumber\\
&\leq&\sum^{\infty}_{n=2}n(|a_n|+|b_n|)|z|^{n-1} \leq \frac{4L_1}{\pi}\sum^{\infty}_{n=2}n|z|^{n-1}
=\frac{4L_1}{\pi}\frac{|z|(2-|z|)}{(1-|z|)^2}.
\eeq
Then by Lemmas \ref{lem-1} and \ref{lem-2}, and \eqref{eq-P}, for $z\in\ID_{r^{\star}}$, we obtain
\beq\label{eq-P1}
|w_z(z)-w_z(0)|&\leq& |c|\, |z|+\left|\frac{\partial}{\partial z}\mathcal{P}_{\gamma_0}(z)-\frac{\partial}{\partial z}\mathcal{P}_{\gamma_0}(0)\right|
+\left|\frac{\partial}{\partial z}\mathcal{G}_1[\gamma](z)-\frac{\partial}{\partial z}\mathcal{G}_1[\gamma](0)\right|\nonumber\\
&&+\left|\frac{\partial}{\partial z}\mathcal{G}_2[g](z)-\frac{\partial}{\partial z}\mathcal{G}_2[g](0)\right|\nonumber\\
&\leq& |z|\left(\frac{4L_1}{\pi}\frac{2-|z|}{(1-|z|)^2}+L_5\right)+4L_3\log\frac{1+|z|}{1-|z|},
\eeq
and
\beq\label{eq-P2}
|w_{\overline{z}}(z)-w_{\overline{z}}(0)|&\leq& |c|\, |z|+\left|\frac{\partial}{\partial {\overline{z}}}\mathcal{P}_{\gamma_0}(z)-\frac{\partial}{\partial {\overline{z}}}\mathcal{P}_{\gamma_0}(0)\right|
+\left|\frac{\partial}{\partial {\overline{z}}}\mathcal{G}_1[\gamma](z)-\frac{\partial}{\partial {\overline{z}}}\mathcal{G}_1[\gamma](0)\right|\nonumber\\
&&+\left|\frac{\partial}{\partial {\overline{z}}}\mathcal{G}_2[g](z)-\frac{\partial}{\partial {\overline{z}}}\mathcal{G}_2[g](0)\right|\nonumber\\
&\leq& |z|\left(\frac{4L_1}{\pi}\frac{2-|z|}{(1-|z|)^2}+L_5\right)+4L_3\log\frac{1+|z|}{1-|z|},
\eeq
where $L_5=|c|+L_2M_1+L_3M_2$.

It follows from Theorem \ref{thm-2}  that
$$1=J_w(0)=\|D_w(0)\|\lambda(D_w(0))\leq \lambda(D_w(0))L_4,
$$
which gives
\beq\label{eq-P3}
\lambda(D_w(0))\geq \frac{1}{L_4}.
\eeq

Now, we are ready to finish the proof of the theorem. First, we demonstrate the univalence of the function $w$ in $\mathbb{D}_{r_{0}}$, where $r_{0}$ is
determined by Eqn. \eqref{fil-1}. For this, let $z_{1}$, $z_{2}$ be two points in $\mathbb{D}_{r_{0}}$ with $z_{1}\neq z_{2}$, and denote the
segment from $z_{1}$ to $z_{2}$ with the endpoints $z_{1}$ and $z_{2}$ by $[z_{1},z_{2}]$. Since
\begin{eqnarray*}
|w(z_{2})-w(z_{1})|&=&\left|\int_{[z_{1},z_{2}]}w_{z}(z)dz+w_{\overline{z}}(z)d\overline{z}\right|
\\\nonumber
&\geq&\left|\int_{[z_{1},z_{2}]}w_{z}(0)dz+w_{\overline{z}}(0)d\overline{z}\right|\\\nonumber
&&
-\left|\int_{[z_{1},z_{2}]}[w_{z}(z)-w_{z}(0)]dz+[w_{\overline{z}}(z)-w_{\overline{z}}(0)]d\overline{z}\right|
,\end{eqnarray*}
we see from \eqref{eq-P1}, \eqref{eq-P2}, \eqref{eq-P3} and Lemma \ref{lem-3} that

\begin{eqnarray*}
|w(z_{2})-w(z_{1})|
&\geq&\lambda(D_{w}(0))\cdot|z_{2}-z_{1}|\\\nonumber
&&-\left|\int_{[z_{1},z_{2}]}[w_{z}(z)-w_{z}(0)]dz+[w_{\overline{z}}(z)-w_{\overline{z}}(0)]d\overline{z}\right|\\\nonumber
&>&\left[\frac{1}{L_4}-2r_0\Big(\frac{4L_1}{\pi}\frac{2-r_0}{(1-r_0)^2}+L_5\Big)-8L_3\log\frac{1+r_0}{1-r_0}\right]|z_{2}-z_{1}|\\\nonumber
&=&0.
\end{eqnarray*}
which implies the univalence of $w$ in $\mathbb{D}_{r_{0}}$.

Next, we prove Theorem \ref{thm-3}\eqref{eq1.10}. For any $\zeta=r_{0}e^{i\theta}\in \partial\mathbb{D}_{r_{0}}$, we obtain that

\begin{eqnarray*}
|w(\zeta)-w(0)|&=&\left|\int_{[0,\zeta]}w_{z}(z)\,dz+w_{\overline{z}}(z)\,d\overline{z}\right|
\\\nonumber
&\geq&\left|\int_{[0,\zeta]}w_{z}(0)dz+w_{\overline{z}}(0)d\overline{z}\right|\\ \nonumber
&&-\left|\int_{[0,\zeta]}[w_{z}(z)-w_{z}(0)]\,dz+[w_{\overline{z}}(z)-w_{\overline{z}}(0)]\,d\overline{z}\right|
\\\nonumber
&\geq&\lambda(D_{w}(0))r_{0}-\frac{8L_1}{\pi}\int_{0}^{r_0}\frac{|z|(2-|z|)}{(1-|z|)^{2}}\,|dz|
-2L_5\int_{0}^{r_0}|z|\,|dz|\\
&&-8L_3\int_{0}^{r_0}\log\frac{1+|z|}{1-|z|}|dz|\hspace{1cm}(\text{by \eqref{eq-P1} and \eqref{eq-P2}})
\\\nonumber
&\geq& \frac{r_{0}}{L_4}-\frac{8L_1}{\pi}\,\frac{r_{0}^2}{1-r_{0}}-L_5r_0^2-8L_3r_0\log\frac{1+r_0}{1-r_0}+\frac{8L_3r_0^2(3-r_0^2)}{(1-r_0^2)^2}
\\\nonumber
&=& \frac{8L_1}{\pi}\,\Big(\frac{r_{0}}{1-r_{0}}\Big)^2+L_5r_0^2+\frac{8L_3r_0^2(3-r_0^2)}{(1-r_0^2)^2} \hspace{1cm}(\text{by \eqref{fil-1}}).
\end{eqnarray*}
Hence $f(\mathbb{D}_{r_{0}})$ contains a univalent disk $\mathbb{D}_{R_{0}}$, where
$$
R_{0}\geq\frac{8L_1}{\pi}\,\Big(\frac{r_{0}}{1-r_{0}}\Big)^2+L_5r_0^2+\frac{8L_3r_0^2(3-r_0^2)}{(1-r_0^2)^2}.
$$
The proof of this theorem is complete.\qed

\subsection*{Acknowledgements}
The research was partly supported by the National Natural Science Foundation of China (No. 11801159).
The  work of the third author is supported by Mathematical Research Impact Centric Support (MATRICS) of
the Department of Science and Technology (DST), India  (MTR/2017/000367).

\subsection*{Conflict of Interests}
The authors declare that there is no conflict of interests regarding the publication of this paper.

\subsection*{Data availability statement}
I confirm that there is no data involved in this manuscript.

\end{document}